\documentclass[12pt]{article}
\usepackage[dvips]{epsfig}
\usepackage{graphics}
\usepackage{latexsym}
\usepackage{verbatim}
\usepackage{amsmath}
\usepackage{amsthm}
\usepackage{amssymb}

\newtheorem{theorem}{Theorem}[section]
\newtheorem{lemma}[theorem]{Lemma}

\def\Remark{\medskip\noindent{\bf Remark: }}
\def\Remarks{\medskip\noindent{\bf Remarks: }}

\newcommand{\ens}[1]{\mathbb{#1}}
\newcommand{\ron}[1]{\mathcal{#1}}

\newcommand{\R}{\mathbb{R}}

\def\cal{\mathcal}

\def\derpar#1#2{\frac{\partial#1}{\partial#2}}
\def\var{\varepsilon}

\def\signcm{\bigskip\bigskip\hspace{80mm}
\vbox{{\sc C. Mouhot\par\vspace{3mm}
UMPA, ENS Lyon\par
46 all\'ee d'Italie\par
69364 Lyon Cedex 07\par
FRANCE\par\vspace{3mm}
e-mail:} cmouhot@umpa.ens-lyon.fr }}

\def\signcb{\bigskip\bigskip\hspace{80mm}
\vbox{{\sc C. Baranger\par\vspace{3mm}
CMLA, ENS Cachan\par
61 avenue du Pr\'esident Wilson\par
94235 Cachan Cedex \par
FRANCE\par\vspace{3mm}
e-mail:} Celine.Baranger@cmla.ens-cachan.fr }}

\begin{document}

\title
{Explicit spectral gap estimates for the linearized Boltzmann and Landau operators with 
hard potentials}

\author{C\'eline Baranger and Cl\'ement Mouhot}

\hyphenation{bounda-ry rea-so-na-ble be-ha-vior pro-per-ties
cha-rac-te-ris-tic}

\maketitle

\begin{abstract} This paper deals with explicit spectral gap estimates for 
the linearized Boltzmann operator with hard potentials (and hard
spheres). We prove that it can be reduced to the Maxwellian case, for
which explicit estimates are already known. Such a method is
constructive, does not rely on Weyl's Theorem and thus does not
require Grad's splitting. The more physical idea of the proof is to
use geometrical properties of the whole collision operator. In a
second part, we use the fact that the Landau operator can be expressed
as the limit of the Boltzmann operator as collisions become grazing in
order to deduce explicit spectral gap estimates for the linearized
Landau operator with hard potentials.
\end{abstract}



\tableofcontents

\section{Introduction}
\setcounter{equation}{0}

This paper is devoted to the study of the spectral properties of the linearized 
Boltzmann and Landau collision operators with hard potentials. In this work we 
shall obtain new quantitative estimates on the spectral gap of these operators. 
Before we explain our methods and results in more details, let us introduce 
the problem in a precise way. The Boltzmann equation describes the behavior of a dilute 
gas when the only interactions taken into account are binary elastic collisions. 
It reads in $\R^N$ ($N \ge 2$) 
 \begin{equation*}
 \derpar{f}{t} + v \cdot \nabla_x f = Q^{{\ron Bo}}(f,f),
 \end{equation*}
where $f(t,x,v)$ stands for the time-dependent distribution function of density of particles 
in the phase space. 
The $N$-dimensional Boltzmann collision operator $Q$ is a quadratic operator, which is local 
in $(t,x)$. The time and position are only parameters and therefore 
shall not be written in the sequel: the estimates proven in this paper are all local in $(t,x)$. 
Thus it acts on $f(v)$ by 
\begin{equation*}
 Q^{{\ron Bo}} (f,f)(v) = \int_{v_* \in \R^N} 
 \int_{\sigma \in \ens{S}^{N-1}}  B(\cos \theta,|v-v_*|) \, 
 \left[ f'_* f' - f_* f \right] \, d\sigma \, dv_*
 \end{equation*} 
where we have used the shorthands $f = f(v)$, $f_* = f(v_*)$, $f ^{'} = f(v')$, 
$f_* ^{'} = f(v_* ^{'})$. The velocities are given by 
 \begin{equation}\label{eq:rel:vit}
 v' = \frac{v+v_*}{2} + \frac{|v-v_*|}{2} \sigma, \qquad
 v'_* = \frac{v+v^*}{2} - \frac{|v-v_*|}{2} \sigma\nonumber.
 \end{equation}
The collision kernel $B$ is a non-negative function which 
only depends on $|v-v_*|$ and $\cos \theta = k \cdot \sigma$ where $k = (v-v_*)/|v-v_*|$. 

Consider the collision operator obtained by the linearization process around the 
Maxwellian {\em global} equilibrium state denoted by $M$
 \begin{equation*}
 L^{{\ron Bo}} h(v)=\int_{\R^N} \int_{\ens{S}^{N-1}} B(\cos \theta,|v-v_*|) \, 
 M(v_{*}) \left[ h^{'}_{*}+h^{'} - h_{*} - h \right] \, d\sigma \, dv_{*},
 \end{equation*}
where $f = M (1+h)$ and $M(v) = e^{-|v|^2}$.
Notice that for the associated linearized equation, the entropy 
is the $L^2(M)$ norm of $h$ and thus by differentiating, 
the entropy production is 
 \begin{equation*}
 <h,L^{{\ron Bo}} h>_{L^2(M)}=-\frac{1}{4}\int_{\R^N}\int_{\R^N}
 \int_{\ens{S}^{N-1}} B(\cos \theta,|v-v_*|) \, 
  \, \left[ h^{'}_{*}+h^{'}-h_{*}-h \right]^2 
 M\, M_{*} \, d\sigma \, dv_{*} \, dv.
 \end{equation*}
This quantity is non-positive. At the level of the 
linearized equation, this corresponds to the first part of Boltzmann's 
$H$-theorem, and it implies that the spectrum of $L^{{\ron Bo}}$ 
in $L^2(M)$ is non-positive. 
Let us denote $D^{{\ron Bo}} (h)=-<h,L^{{\ron Bo}} h>_{L^2(M)}$. We shall 
call this quantity ``(linearized) entropy dissipation functional'' by analogy 
with the nonlinear case.

In the case of long-distance interaction, the collisions occur 
mostly for very small deviation angle $\theta$. In the case of the  
Coulomb potential, for which the Boltzmann collision operator is 
meaningless (see~\cite[Annex I, Appendix A]{Vill:habil:00}), one has to 
replace it by the Landau collision operator
 \begin{equation*}
 Q^{{\ron La}} (f,f)(v) = \nabla _v \cdot \left( \int_{v_* \in \R^N} 
 {\bf a}(v-v_*) \left[ f_* \left( \nabla f \right) 
 - f \left( \nabla f \right)_* \right] \, dv_* \right),
 \end{equation*} 
with ${\bf a} (z) = |z|^2 \, \Phi(z) \, \Pi_{z^\bot}$, where $\Pi_{z^\bot}$ is 
the orthogonal projection onto $z^\bot$, i.e 
 \begin{equation*}
 \left( \Pi_{z^\bot} \right) _{i,j} = \delta_{i,j} - \frac{z_i z_j}{|z|^2}.
 \end{equation*}

This operator is used for instance in models of plasma in the case of a Coulomb 
potential, i.e a gas of (partially or 
totally) ionized particles (for more details see~\cite{Vill:hand} and the references 
therein). Applying the same linearization process than for the Boltzmann operator (around 
the same global equilibrium $M$), we define the linearized Landau operator
 \begin{equation*}
 L^{{\ron La}} h(v)=M(v) ^{-1} \, \nabla_v \cdot \left( \int_{v_* \in \R^N} 
 {\bf a}(v-v_*) \left[ \left( \nabla h \right)  - \left( \nabla h \right)_* \right] 
 M M_* \, dv_* \right),
 \end{equation*}
and the (linearized) Landau entropy dissipation functional 
 \begin{eqnarray*}\nonumber
 D^{{\ron La}} (h)& =& - <h,L^{{\ron La}}h>_{L^2(M)} \\
 &=&\frac{1}{2}\int_{\R^N}\int_{\R^N}\Phi(v-v_*)|v-v_*|^2
 \left\|\Pi_{(v-v_*)^\bot} \left[ \left( \nabla h \right) - \left( \nabla h \right)_*\right] \right\|^2 
 M\, M_{*} \, dv_{*}\, dv
 \end{eqnarray*}
which is also non-positive. It implies that the spectrum of $L^{{\ron La}}$ in $L^2(M)$ is non-positive.

\medskip
Let us now write down our assumptions for the collision kernel $B$:
 \begin{itemize}
 \item $B$ is a tensorial product
  \begin{equation}\label{eq:hyptens}
  B = b(\cos \theta) \, \Phi(|v-v_*|),
  \end{equation}
 where $\Phi$ and $b$ are non-negative functions (this is the case for instance for 
 collision kernels deriving from interaction potentials behaving like inverse-power laws).
 \item The kinetic part $\Phi$ is bounded from below at infinity, i.e 
  \begin{equation}\label{eq:hypPhi}
  \exists \, R \ge 0, \ c_\Phi >0 \ \  | \ \ \forall \, r \ge R, \ \Phi(r) \ge c_\Phi.
  \end{equation}
 This assumption holds for hard potentials (and hard spheres).
 \item The angular part $b$ satisfies
 \begin{equation}\label{eq:hypb}
  c_b = \inf_{\sigma_1, \sigma_2 \in \ens{S}^{N-1}} 
  \int_{\sigma_3 \in \ens{S}^{N-1}} \min \{  b(\sigma_1 \cdot \sigma_3), 
  b(\sigma_2 \cdot \sigma_3) \} \, d\sigma_3 > 0. 
  \end{equation}
 This covers all the physical cases.
 \end{itemize}
\Remarks 1. Notice that there is no $b$ left in $Q^{{\ron La}}$ and $L^{{\ron La}}$ but 
the function $\Phi$ is definitely the same in both Landau and Boltzmann operators. 
Therefore the assumptions on the Landau operator reduce to~\eqref{eq:hypPhi}. Thus  
we deal with the so-called ``hard potentials'' case for the Landau operator, which 
excludes the Coulomb potential.
\smallskip

2. The assumption that $B$ is a tensorial product is made for a sake of 
simplicity. Indeed, one could easily adapt 
the proofs in section~\ref{sec:Bolt} to relax this assumption. 
The price to pay would be a more technical 
condition on the collision kernel $B$. 
\medskip

The spectral properties of the linearized Boltzmann and Landau
operators have been extensively studied. In particular, there are of
crucial interest for perturbative approach issues. For instance, the
convergence to equilibrium has been studied in this context, as well
as the hydrodynamical limit (see~\cite{ElPi:LBE:75}).

\begin{figure}[h]
 \epsfysize=4cm
 $$\epsfbox{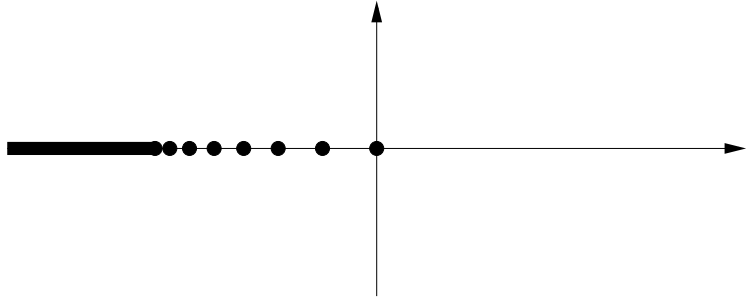}$$
 \caption{Spectrum of the collision operator for strictly hard potential 
 with angular cutoff}\label{fig:spec:maxw}
 \end{figure}

\begin{figure}[h]
 \epsfysize=4cm
 $$\epsfbox{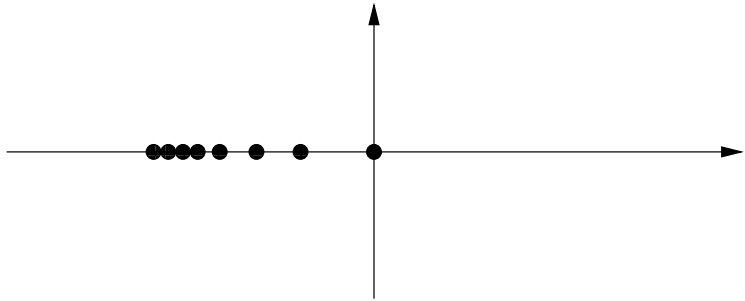}$$
 \caption{Spectrum of the collision operator for Maxwell's molecules 
 with angular cutoff}\label{fig:spec:hard}
 \end{figure}

On the one hand, for hard potentials, the existence of a spectral gap as soon 
as the kinetic part of the collision kernel is bounded from below at infinity is 
a classical result, which can be traced back unto Grad himself. The only method 
was up to now to work under the assumption of Grad's angular cutoff, and to apply Weyl's 
Theorem to $L^\cal{B}$, written as a compact perturbation of a multiplication operator 
(a very clear presentation of this proof can be found in~\cite{CIP:94}). 
The picture of the spectrum obtained for the operator (under Grad's 
cutoff assumption) is described by figures~\ref{fig:spec:maxw} and~\ref{fig:spec:hard} 
(see~\cite{Cerc:EB:88}). 

A similar method has been applied to the Landau linearized
operator with hard potential in~\cite{DeLe:LLE:97}.

On the other hand, for the particular case of Maxwellian molecules (for 
$L^\cal{B}$), a complete and explicit 
diagonalisation has been obtained first by symmetry arguments in~\cite{WCUh:LBE:70}, and 
then by Fourier methods in~\cite{Boby:maxw:88}. 
The spectral gap for the ``over-Maxwellian'' collision kernel of the Landau 
linearized operator (i.e collision kernels which 
are bounded from below by one for Maxwellian molecules) can be derived from results
in~\cite{DeVi:LII:00}, by a linearization process.   
Notice also that in the case of the so-called Kac's equation, an explicit entropy
production estimate, based on a cancellation method, was given in~\cite{Desv:cveq:98}; 
this method can be linearized in order to give explicit spectral gap estimates for ``over-quadratic'' 
linearized Kac's operator (for which the physical meaning is not clear!). Nevertheless we 
did not manage to adapt this strategy to the Boltzmann operator with hard potentials. 
Notice however that Wennberg~\cite{Wenn:entrop} gave an extension of the very first entropy estimates of 
Desvillettes~\cite{Desv:cveq:89} to allow for hard and soft potentials. His 
idea has some similarities with ours: to avoid the region in $\R^N \times \R^N$ 
where $\Phi(|v-v_*|)$ is small. 

A specific study of the spectral properties of the linearized operator 
was made for {\em non-cutoff} hard potentials in~\cite{Pao:TSnoncutoff:74}. 
Nevertheless this article was critically reviewed some years later in~\cite{Klau:TSnoncutoff:77}. 
Also another specific study for ``radial cutoff potentials'' was done in~\cite{CzPa:radcutoff:80}. 

\begin{figure}[h]
 \epsfysize=4cm
 $$\epsfbox{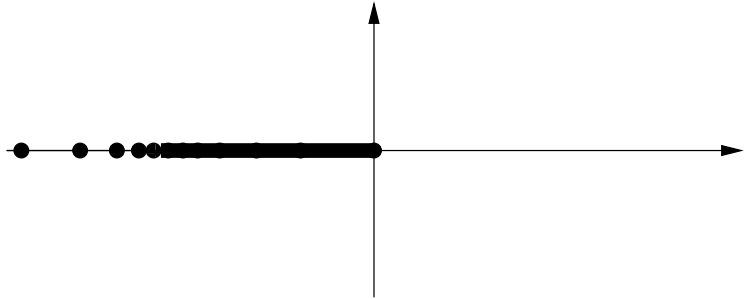}$$
 \caption{Spectrum of the collision operator for soft potentials 
 with angular cutoff}\label{fig:spec:soft}
 \end{figure}

Finally notice that it is proved in~\cite{Cafl:LBE1:80} that the Boltzmann linearized 
operator with soft potential has no spectral gap. The resulting spectrum is described in 
figure~\ref{fig:spec:soft} (see~\cite{Cerc:EB:88}). 

But if one allows a loss on the algebraic weight of the norm, 
it was proved in~\cite{GolPo:cras:86} a ``degenerated spectral gap'' result 
of the form 
 \[ \| L^{{\ron Bo}} h \|_{L^2(M)} \ge C \, \| h \|_{L^2 _\gamma(M)}  \ \ 
 \forall \, h \bot \left\{ 1 ; v \, ; |v|^2 \right\}, \] 
where $\gamma <0$ is the power of the kernel $\Phi$. It is based on 
inequalities proven in~\cite{Cafl:LBE1:80} and Weyl's Theorem. 

\bigskip

However, the perturbative method has drawbacks, all coming from the 
fact that it does not rely on a physical argument.
First it is not explicit, that is the width of the spectral gap 
is not known, which is problematic when one wants to
obtain quantitative estimates of convergence to equilibrium.
Secondly it gives no information about how this spectral gap is 
sensitive to the perturbation of the collision kernel. Finally approaches 
based on Weyl's Theorem rely strongly on Grad's cutoff assumption via 
``Grad's splitting'', which means to deal separately with the gain and the 
loss part of the collision operator.

Our method is geometrical and based on a physical argument. It gives explicit
estimates and deals with the whole operator, with or without angular cutoff.  
Up to our knowledge, as far as spectral gaps are considered, it covers 
all the results of the above-mentioned articles dealing with hard potentials, 
with or without angular cutoff. \\
We think likely that this geometrical method could also be adapted to give 
explicit versions of ``degenerated spectral gap'' results in the case 
of soft potentials, even if up to now we did not manage to do it.

\smallskip 

We now state our main theorems:
 \begin{theorem}[The Boltzmann linearized operator] \label{theo:Bolt}
 Under the assumptions~\eqref{eq:hyptens}, \eqref{eq:hypPhi}, \eqref{eq:hypb}, 
 the Boltzmann entropy dissipation functional $D^{{\ron Bo}}$ with $B = \Phi \, b$ 
 satisfies, for all $h \in L^2(M)$
  \begin{equation}\label{eq:reducB}
  D^{{\ron Bo}} (h) \ge C^{{\ron Bo}} _{\Phi,b} \, D^{{\ron Bo}} _0 (h),
  \end{equation}
 where $D^{{\ron Bo}} _0 (h)$ stands for the entropy dissipation functional with 
 $B_0 \equiv 1$ and
  \begin{equation*}
  C^{{\ron Bo}} _{\Phi,b} = \left( \frac{c_\Phi \, c_b \, e^{-4 R^2}}{32 \left|\ens{S}^{N-1}\right|} \right)  
  \end{equation*}
 with $R$, $c_\Phi$, $c_b$ being defined in~\eqref{eq:hypPhi}, \eqref{eq:hypb}.
 
 As a consequence we deduce quantitative estimates on the spectral gap of the 
 linearized Boltzmann operator, namely for all $h \in L^2 (M)$ 
 orthogonal in $L^2(M)$ to $1$, $v$ and $|v|^2$, we have 
  \begin{equation}\label{eq:spB}
  D^{{\ron Bo}} (h) \ge C^{{\ron Bo}} _{\Phi. b} \, |\lambda^{{\ron Bo}} _0| \, \| h \|^2 _{L^2(M)}.
  \end{equation}
 Here $\lambda^{{\ron Bo}} _0$ is the first non-zero eigenvalue of the 
 linearized Boltzmann operator with $B_0 \equiv 1$
 (that is, for Maxwellian molecules with no angular dependence, sometimes 
 called pseudo-Maxwellian molecules) which equals in dimension $3$ 
 (see~\cite{Boby:maxw:88})
  \begin{equation*} 
  \lambda^{{\ron Bo}} _0 = - \pi \, \int_0 ^\pi \sin ^3 \theta \, d\theta = - \frac{4 \pi}{3}. 
  \end{equation*}
 \end{theorem}

\Remark As an application of this theorem, let us give explicit  
formulas for the spectral gap $S_\gamma ^{{\ron Bo}}$ of the 
Boltzmann linearized operator with $b \ge 1$ 
and $\Phi(z) = |z|^\gamma$, $\gamma >0$, in dimension $3$. 
Then $c_b \ge |\ens{S}^2|$ 
and for any given $R$ we can take $c_\Phi = R^\gamma$. Thus we get 
 \[ S_\gamma ^{{\ron Bo}} \ge \left( \frac{R^\gamma \, e^{-4 R^2}}{32} \right) 
 \,  \frac{4 \pi}{3}  \] 
for any $R >0$. An easy computation leads to the lower bound 
 \[ S_\gamma ^{{\ron Bo}} \ge \frac{\pi \, (\gamma/8)^{\gamma/2} \, e^{-\gamma/2}}{24} \]
by optimizing the free parameter $R$.
\medskip

 \begin{theorem}[The Landau linearized operator] \label{theo:Land}
 Under assumptions~\eqref{eq:hypPhi}, the Landau entropy dissipation functional $D^{{\ron La}}$ 
 with collision kernel $\Phi$ satisfies, for all $h \in L^2(M)$
  \begin{equation}\label{eq:reduc:Landau}
  D^{{\ron La}}(h) \ge C^{{\ron La}} _\Phi \, D^{{\ron La}} _0 (h)
  \end{equation}   
 where $D^{{\ron La}} _0 (h)$ stands for the Landau entropy dissipation functional 
 with $\Phi_0 \equiv 1$  and 
  \[ C^{{\ron La}} _\Phi = \left( \frac{c_\Phi \, \beta_R}{8 \, \alpha_N} \right) \] 
 with 
  \[ \alpha_N = \int_{\R^{N-1}} e^{-|V|^2} \, dV, \ \ \ \ 
     \beta_R = \int_{\big\{ V \in \R^{N-1} \ | \ |V| \ge 2R \big\}} 
               e^{-|V|^2} \, dV. \]
 
 As a consequence we deduce quantitaves estimates on the spectral gap of the 
 linearized Landau operator, namely for all $h \in L^2 (M)$ 
 orthogonal in $L^2(M)$ to $1$, $v$ and $|v|^2$, we have 
  \begin{equation}\label{eq:tsL}
  D^{{\ron La}} (h) \ge C^{{\ron La}} _{\Phi} \, |\lambda^{{\ron La}} _0| \, \| h \|^2 _{L^2(M)}.
  \end{equation}
 Here $\lambda^{{\ron La}} _0$ is the first non-zero eigenvalue of the 
 linearized Landau operator with $\Phi_0 \equiv 1$
 (that is, for Maxwellian molecules).

 Moreover in dimension $3$, by grazing collisions limit, we can estimate 
 $\lambda^{{\ron La}} _0$ thanks to the explicit formula on the spectral gap 
 of the Boltzmann linearized operator for Maxwellian molecules
  \begin{equation}\label{eq:vpL}
  |\lambda^{{\ron La}} _0| \ge 2\, \pi.
  \end{equation}
 \end{theorem}

\Remarks
1. As for the Boltzmann linearized operator, we can deduce from this theorem 
an explicit formula for a lower bound on the spectral gap $S_\gamma ^{{\ron La}}$ 
for the Landau linearized operator with hard potentials 
$\Phi(z) = |z|^\gamma$, $\gamma >0$, in dimension $3$. 
We get 
 \[ S_\gamma ^{{\ron La}} \ge \left( \frac{R^\gamma \, e^{-4 R^2}}{8} \right) 
 \,  2 \pi \] 
for any $R >0$. An easy computation leads to the lower bound 
 \[ S_\gamma ^{{\ron La}} \ge \frac{\pi \, (\gamma/8)^{\gamma/2} \, e^{-\gamma/2}}{4} \]
by optimizing the free parameter $R$. 
\smallskip

2. The modulus of the first non-zero eigenvalue of the Landau linearized operator 
for Maxwellian molecules is estimated here by grazing collisions limit. 
Other methods would have been the linearization of entropy estimates in~\cite{DeVi:LII:00}, or  
to use the decomposition (established in~\cite{Vill:LandMax:98})
of the Landau operator for Maxwellian molecules 
into a Fokker-Planck part (for which spectral gap is already known) and 
a spherical diffusion process, which can only increases the spectral gap; and 
then to linearize the estimate thus obtained.
\smallskip

3. More generally, it is likely that an explicit spectral gap for the 
Landau linearized operator with hard potentials could be directly  
computed by existing methods even if up to our knowledge this is the first explicit formula.  
But Theorem~\ref{theo:Land} is stronger: it says that the property proved on 
the Boltzmann operator with hard potentials, namely ``cancellations for small relative velocities 
can be neglected as far as entropy production is concerned'', remains true for the Landau 
linearized operator with hard potentials. 
\medskip

The idea of the proof is to reduce the case of hard potentials 
(in the generalized sense~\eqref{eq:hypPhi}) 
to the Maxwellian case. 
The difficulty is to deal with the cancellations of the kinetic collision kernel $\Phi$ 
on the diagonal $v = v_*$. 

The starting point is the following inequality which is a corollary
 of~\cite[Theorem 2.4]{CMCV:gran:01}
 \begin{multline}\label{eq:CMCV}
 \int_{\R^N} \int_{\R^N} | \xi (x) - \xi(y)|^2 \, |x-y|^\gamma \, M(x) \, M(y) \, dx \, dy \\
 \ge K_\gamma \,  \int_{\R^N} \int_{\R^N} | \xi (x) - \xi(y)|^2  \, M(x) \, M(y) \, dx \, dy
 \end{multline}
for $\gamma \ge 0$, $\xi$ some function, and 
 \begin{equation*}
 K_\gamma = \frac{1}{4 \int_{\R^N} M} \ \inf_{x,y \in \R^N} \int_{\R^N} 
 \min \left\{ |x-z|^\gamma , |z-y|^\gamma \right\} \, M(z) \, dz .
 \end{equation*}

It was first suggested by Villani~\cite[Chap. 5, section 1.4]{Vill:hand}, in the 
context of the study of entropy-entropy dissipation inequalities for the 
Landau equation with hard potentials, that this inequality could allow to 
prove that hard potentials reduce to the Maxwellian case as far as 
convergence to equilibrium is concerned. 

The proof of~\eqref{eq:CMCV} relies strongly on the existence 
of a ``triangular inequality'' for some function $F(x,y)$ integrated: 
in~\eqref{eq:CMCV}, the function $F$ is simply $|\xi(x) - \xi(y)|^2$ 
which satisfies 
 \[ F(x,y) \le 2 F(x,z) + 2 F(z,y). \]
The main difficulty is hence to obtain such a ``triangular inequality'' adapted to our case for the 
Boltzmann linearized operator. It will be discussed in details in section~\ref{sec:Bolt} together 
with the proof of Theorem~\ref{theo:Bolt}. Section~\ref{sec:Landau} will be devoted to the Landau linearized 
operator: using results of section~\ref{sec:Bolt}, we will prove Theorem~\ref{theo:Land} 
thanks to a grazing collision limit.
 
\section{The Boltzmann linearized operator} \label{sec:Bolt}
\setcounter{equation}{0}
In this section, we present the proof of inequality~\eqref{eq:reducB} in Theorem~\ref{theo:Bolt}. 
In order to ``avoid'' the diagonal $v \sim v_*$ where $\Phi$ is not uniformly bounded from below,
we use the following argument: 
perfoming a collision with small relative velocity (i.e. for a small $|v-v_*|$) is the same than 
perfoming two collisions with great relative velocity, provided that the pre- and 
post-collisionnal velocities are the same. One could summarize the situation in this 
way: when a collision with small relative velocity occurs, at the same time, two collisions with great 
relative velocity occur, which give the same pre- and post-collisionnal velocities, 
and which produce {\em at least the same amount of entropy}.

Before proving~\eqref{eq:reducB}, let us begin with a preliminary lemma dealing with the
angular part of the collision kernel. This lemma is based on the same geometrical idea as the one we shall 
use for the treatment of the cancellations of $\Phi$: the introduction of some well-chosen 
intermediate collision.     
This first step is made for the sake of simplicity: we show that 
in the sequel of this section one can set $b \equiv 1$ without restriction. It 
makes the proof clearer, and simplifies somehow the constants. 

Let us denote from now on 
 \[ k (v,v_*,v',v' _*) = \big[ h(v) + h(v_*) - h(v') - h(v' _*) \big]^2. \]
 
 \begin{lemma}[Homogeneization of the angular collision kernel $b$]\label{lem:b}
 Under the assumptions~\eqref{eq:hyptens}, \eqref{eq:hypPhi}, \eqref{eq:hypb}, 
 for all $h \in L^2 (M)$,
  \begin{equation}\label{eq:ineg1}
  D^{{\ron Bo}} (h) \ge \frac{c_b}{4 \left|\ens{S}^{N-1}\right|} D^{{\ron Bo}} _1 (h)
  \end{equation}
 where $D^{{\ron Bo}} _1$ denotes the entropy dissipation functional with $B = \Phi (|v-v_*|)$ 
 instead of $B = \Phi (|v-v_*|) \, b (\theta)$.
 \end{lemma}

\Remark This lemma allows to bound from below the entropy dissipation functional by 
one with an ``uniform angular collision kernel'', i.e a constant $c_b$, 
even when $b$ is not bounded from below by a positive number uniformly on the 
sphere. Notice for instance that the condition $c_b>0$ is satisfied for $b$ 
having only finite number of $0$. 

\begin{proof}[Proof of Lemma~\ref{lem:b}]
First, we write down an appropriate representation of the operator. 
The functional $D^{{\ron Bo}}$ reads in ``$\sigma$-representation''
 \begin{equation*}
 D^{{\ron Bo}} (h) = \frac{1}{4}\,\int_{\R^N}\int_{\R^N}
 \int_{\ens{S}^{N-1}} \Phi(|v-v_*|) \, b \left(\frac{v-v_*}{|v-v_*|} \cdot \sigma_1 \right) \, 
 M\, M_{*} \, k(v,v_*,v',v' _*)
 \, d\sigma_1 \, dv_{*}\, dv.
 \end{equation*}
(for the classical representations of the Boltzmann operator we refer to~\cite{Vill:hand}). 
Then keeping $\sigma_1$ fixed we do the change of variable 
$(v,v_*) \rightarrow (\frac{v+v_*}{2},\frac{v-v_*}{2})$, whose jacobian is $(-1/2)^N$. Let us 
denote $\Omega = \frac{v + v_*}{2}$ and 
$\Omega' = \frac{v -v_*}{2}$. We obtain 
 \begin{eqnarray*} 
 D^{{\ron Bo}} (h) &=& \frac{2^N}{4}\, \int_{\Omega \in \R^N}\int_{\Omega' \in \R^N}
 \int_{\ens{S}^{N-1}} \Phi(2|\Omega'|) \, b\left(\frac{\Omega'}{|\Omega'|}\cdot \sigma_1\right) \\ 
 && k\big( \Omega+\Omega', \Omega-\Omega', 
 \Omega+|\Omega'|\sigma_1, \Omega- |\Omega'| \sigma_1 \big) 
 e^{-2|\Omega|^2-2|\Omega'|^2} \, d\sigma_1 \, d\Omega \, d\Omega'
 \end{eqnarray*}
(recall that $|\Omega|^2+|\Omega'|^2 = \frac{|v|^2 + |v_*|^2}2$).

We now write $\Omega'$ in spherical coordinates $\Omega' = r \, \sigma_2$, 
the other variables being kept fixed, and use Fubini's Theorem
 \begin{eqnarray*} \nonumber
 D^{{\ron Bo}} (h) &=&  \frac{2^N}{4}\, \int_{\Omega \in \R^N} \int_{r \in \R_+}
 r^{N-1} \, \Phi(2r) \, e^{-2|\Omega|^2-2r^2} \int_{\sigma_1 \in \ens{S}^{N-1}} \int_{\sigma_2 \in \ens{S}^{N-1}}
  b(\sigma_1 \cdot \sigma_2) \\ 
 && k \big( \Omega+r\sigma_2, \Omega-r\sigma_2, 
 \Omega+r\sigma_1, \Omega- r\sigma_1 \big) 
  \, d\sigma_1 \, d\sigma_2 \, dr \, d\Omega.
 \end{eqnarray*}
Now we apply a geometrical idea that we shall also use below in the treatment of cancellations 
of $\Phi$: namely we add a third artificial variable. 
Let us thus introduce two collisions points $u$ and $u_*$ on the sphere 
of center $\Omega$ and radius $r$ (see figure~\ref{fig:ang}) and replace
the collision ``$(v,v_*)$ gives $(v',v'_*)$'' by the two collisions
``$(v,v_*)$  gives $(u,u_*)$'' and  ``$(u,u_*)$ gives $(v',v'_*)$''.
 \begin{figure}[h]
 \epsfysize=6cm
 $$\epsfbox{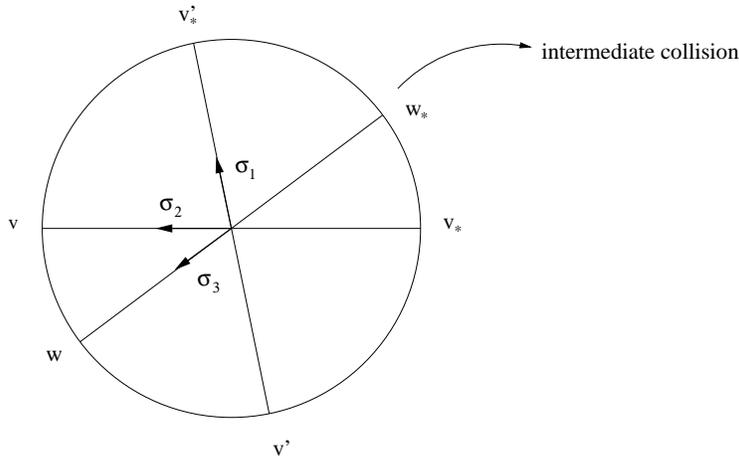}$$
 \caption{Introduction of an intermediate collision}\label{fig:ang}
 \end{figure}

Then, we shall use the following ``triangular'' inequality on the collision points:
 \begin{eqnarray}\label{eq:triang1}
 \big[ \left( h(v) + h(v_*) \right) - \left( h(v') + h(v' _*) \right) \big]^2 
 &\le& 2 \big[ \left( h(v) + h(v_*) \right) - \left( h(u) + h(u_*) \right) \big]^2 \\
 && + 2 \big[ \left( h(u) + h(u_*) \right) - \left( h(v') + h(v'_*) \right) \big]^2 \nonumber. 
 \end{eqnarray}
So let us add a third ``blind'' variable $\sigma_3$ on the sphere
 \begin{eqnarray*} \nonumber
 D^{{\ron Bo}} (h) &=&  \frac{2^N}{4\,\left|\ens{S}^{N-1}\right|} \, 
 \int_{\Omega \in \R^N} \int_{r \in \R_+}
 r^{N-1} \, \Phi(2r) \, e^{-2|\Omega|^2-2r^2} \int_{\sigma_1 \in \ens{S}^{N-1}} 
 \int_{\sigma_2 \in \ens{S}^{N-1}} \int_{\sigma_3 \in \ens{S}^{N-1}} \\  
 && b(\sigma_1 \cdot \sigma_2) 
 \, k \big( \Omega+r\sigma_2, \Omega-r\sigma_2, 
 \Omega+r\sigma_1, \Omega- r\sigma_1 \big) 
  \, d\sigma_1 \, d\sigma_2 \, d\sigma_3 \, dr \, d\Omega.
 \end{eqnarray*}
As variables $\sigma_1$, $\sigma_2$ and $\sigma_3$ are equivalent, one can change the 
``blind'' variable into either $\sigma_1$ or $\sigma_2$ and compute the mean to get
 \begin{eqnarray} \nonumber
 D^{{\ron Bo}} (h) &=&  \frac{2^N}{4\,\left|\ens{S}^{N-1}\right|} \, 
 \int_{\Omega \in \R^N} \int_{r \in \R_+}
 r^{N-1} \, \Phi(2r) \, e^{-2|\Omega|^2-2r^2} 
 \int_{\sigma_1 \in \ens{S}^{N-1}} \int_{\sigma_2 \in \ens{S}^{N-1}}
 \int_{\sigma_3 \in \ens{S}^{N-1}} \nonumber\\ 
 && \frac{1}{2}\,\Big[ b(\sigma_1 \cdot \sigma_3) \,
 k \big( \Omega+r\sigma_3, \Omega-r\sigma_3, 
 \Omega+r\sigma_1, \Omega- r\sigma_1 \big) \nonumber\\
 && + b(\sigma_2 \cdot \sigma_3) \, k \big( \Omega+r\sigma_2, \Omega-r\sigma_2,
 \Omega+r\sigma_3, \Omega- r\sigma_3 \big)\Big]
  \, d\sigma_1 \, d\sigma_2 \, d\sigma_3 \, dr \, d\Omega\nonumber,
 \end{eqnarray} 
which yields
 \begin{eqnarray*}
 D^{{\ron Bo}} (h) &\ge&  \frac{2^N}{4\,\left|\ens{S}^{N-1}\right|} \, 
 \int_{\Omega \in \R^N} \int_{r \in \R_+}
 \int_{\sigma_1 \in \ens{S}^{N-1}} \int_{\sigma_2 \in \ens{S}^{N-1}}
 \int_{\sigma_3 \in \ens{S}^{N-1}} \Phi(2r) \\ 
 && \frac{1}{2}\,\min\{  b(\sigma_1 \cdot \sigma_3), 
 b(\sigma_2 \cdot \sigma_3) \}\,\Big[ k \big( \Omega+r\sigma_3, \Omega-r\sigma_3, 
 \Omega+r\sigma_1, \Omega- r\sigma_1 \big) \\
 &&+ k \big( \Omega+r\sigma_2, \Omega-r\sigma_2,
 \Omega+r\sigma_3, \Omega- r\sigma_3 \big) \Big]
 e^{-2|\Omega|^2-2r^2} \, d\sigma_1 \, d\sigma_2 \, d\sigma_3 \, dr \, d\Omega,
 \end{eqnarray*} 
The triangular inequality needed on $k$ is 
 \begin{eqnarray} \nonumber
 k \big(\Omega+r\sigma_2, \Omega-r\sigma_2, 
 \Omega+r\sigma_1, \Omega- r\sigma_1 \big) &\le& 
 2 \, k \big(\Omega+r\sigma_3, \Omega-r\sigma_3, 
 \Omega+r\sigma_1, \Omega- r\sigma_1 \big) \\ \nonumber
 && + 2 \, k \big( \Omega+r\sigma_2, \Omega-r\sigma_2, 
 \Omega+r\sigma_3, \Omega- r\sigma_3 \big)
 \end{eqnarray}
and follows from~\eqref{eq:triang1}. Thus if one sets
 \[ c_b = \inf_{\sigma_1, \sigma_2 \in \ens{S}^{N-1}} 
 \int_{\sigma_3 \in \ens{S}^{N-1}} \min \{  b(\sigma_1 \cdot \sigma_3), 
 b(\sigma_2 \cdot \sigma_3) \} \, d\sigma_3, \] 
one obtains (going back to the classical representation)
 \begin{eqnarray*} 
 D^{{\ron Bo}} (h) &\ge& \frac{c_b}{4 \left|\ens{S}^{N-1}\right|}\frac{1}{4} 
 \int_{\ens{S}^{N-1}} \int_{\R^N}\int_{\R^N} 
 \Phi(|v-v_*|) \,  M\, M_{*} \, k(v,v_*,v',v' _*) 
 \, d\sigma  \, dv_{*}dv \\ 
 &\ge& \frac{c_b}{4 \left|\ens{S}^{N-1}\right|} D^{{\ron Bo}} _1 (h)
 \end{eqnarray*}
which concludes the proof.
\end{proof}
 \begin{lemma}[Treatment of the cancellations of $\Phi$]\label{lem:phi}
 Under the assumptions~\eqref{eq:hypPhi} on $\Phi$,  
 for all $h \in L^2(M)$ 
  \begin{equation}\label{eq:ineg2} 
  D^{{\ron Bo}} _1 (h) \ge \left( \frac{c_\Phi \, e^{-4 R^2}}{8} \right) D^{{\ron Bo}} _0 (h)
  \end{equation}
 where $D^{{\ron Bo}} _1$ is the entropy dissipation functional with $B = \Phi (|v-v_*|)$ 
 and $D^{{\ron Bo}} _0$ is the entropy dissipation functional with $B=1$.
 \end{lemma}
\begin{proof}[Proof of Lemma~\ref{lem:phi}]
We assume here that $b \equiv 1$. Lemma~\ref{lem:b} indeed shows that this is no 
restriction modulo a factor $c_b/(4 \left|\ens{S}^{N-1}\right|)$. 
Let us consider the so-called ``$\omega$-representation'' (see~\cite{Vill:hand} again): 
the vector $\sigma$ 
integrated on the sphere becomes $\omega = \frac{v'-v}{|v'-v|}$ and the 
change of variable changes the angular kernel into 
 \begin{equation*} 
 \tilde{b} (\theta) = 2^{N-1} \sin^{N-2} \left( \frac{\theta}{2} \right).
 \end{equation*} 
where $\cos \theta = 2 (k \cdot \omega)^2 -1$ with $k = (v-v_*)/|v-v_*|$.\\
The operator $D^{{\ron Bo}} _1 (h)$ thus becomes
 \begin{equation*} 
 D^{{\ron Bo}} _1 (h)=\frac{1}{4}\, \int_{\R^N}\int_{\R^N}
 \int_{\ens{S}^{N-1}} \Phi(|v-v_*|) \, \tilde{b}(\theta) \, 
 M\, M_{*} \, k(v,v_*,v',v' _*)
 \, d\omega \, dv \, dv_*.
 \end{equation*}
where the velocities $v',v' _*$ are given by
 \[ v' = v - ( v -v_*,\omega) \omega, \ \ \ v' _* = v_* + ( v -v_*,\omega) \omega. \]
Then keeping $\omega$ fixed we do the following change of variable
 \[ v = r_1 \omega + V_1, \ \ \ \ v_* = r_2 \omega + V_2 \] 
with $V_1, V_2 \in \omega^\bot$. 
The Jacobian of the change of variable is $1$ since the decompositions 
are orthogonal.
Finally we obtain the following representation
 \begin{eqnarray*} 
 D^{{\ron Bo}} _1 (h) &=&\frac{1}{4} \, \int_{\ens{S}^{N-1}} \int_{V_1 \in \omega^\bot} 
 \int_{V_2 \in \omega^\bot} e^{-|V_1|^2-|V_2|^2} \int_{r_1 \in \R} \int_{r_2 \in \R}
 e^{-r_1 ^2-r_2 ^2} \, \Phi \left(\sqrt{|r_2 - r_1|^2 + |V_2 - V_1|^2} \right) \\ \nonumber
 &&  \, \tilde{b}(\theta) \, 
 k(r_1 \omega + V_1, r_2 \omega + V_2, r_2 \omega + V_1, r_1 \omega + V_2) \,
 dr_1 \, dr_2 \, dV_2 \, dV_1 \, d\omega. 
 \end{eqnarray*}
Assume that $\Phi$ is non-decreasing. 
This is no restriction since $\Phi \ge \tilde{\Phi}$, with 
 \begin{equation*}
 \tilde{\Phi} (r) = \inf_{r' \ge r} \Phi(r'),
 \end{equation*}
and $\tilde{\Phi}$ satisfies assumption~\eqref{eq:hypPhi} 
with the same constant as $\Phi$. This 
monotonicity yields
 \begin{eqnarray*} 
 D^{{\ron Bo}} _1 (h) &\ge&\frac{1}{4} \, \int_{\ens{S}^{N-1}} \int_{V_1 \in \omega^\bot} 
 \int_{V_2 \in \omega^\bot} e^{-|V_1|^2-|V_2|^2} \int_{r_1 \in \R} \int_{r_2 \in \R}
 e^{-r_1 ^2-r_2 ^2} \, \Phi \left(|r_2 - r_1|\right) \\ \nonumber
 &&  \, \tilde{b}(\theta) \, 
 k(r_1 \omega + V_1, r_2 \omega + V_2, r_2 \omega + V_1, r_1 \omega + V_2) \,
 dr_1 \, dr_2 \, dV_2 \, dV_1 \, d\omega. 
 \end{eqnarray*}

We now introduce two collision points $u$ and $u_*$ (see figure~\ref{fig:2D}) and replace
the collision ``$(v,v_*)$ gives $(v',v'_*)$'' by the two collisions
``$(v,u_*)$  gives $(v'_*,u)$'' and  ``$(u,v_*)$ gives $(u_*,v')$''.
 \begin{figure}[h]
 \epsfysize=6cm
 $$\epsfbox{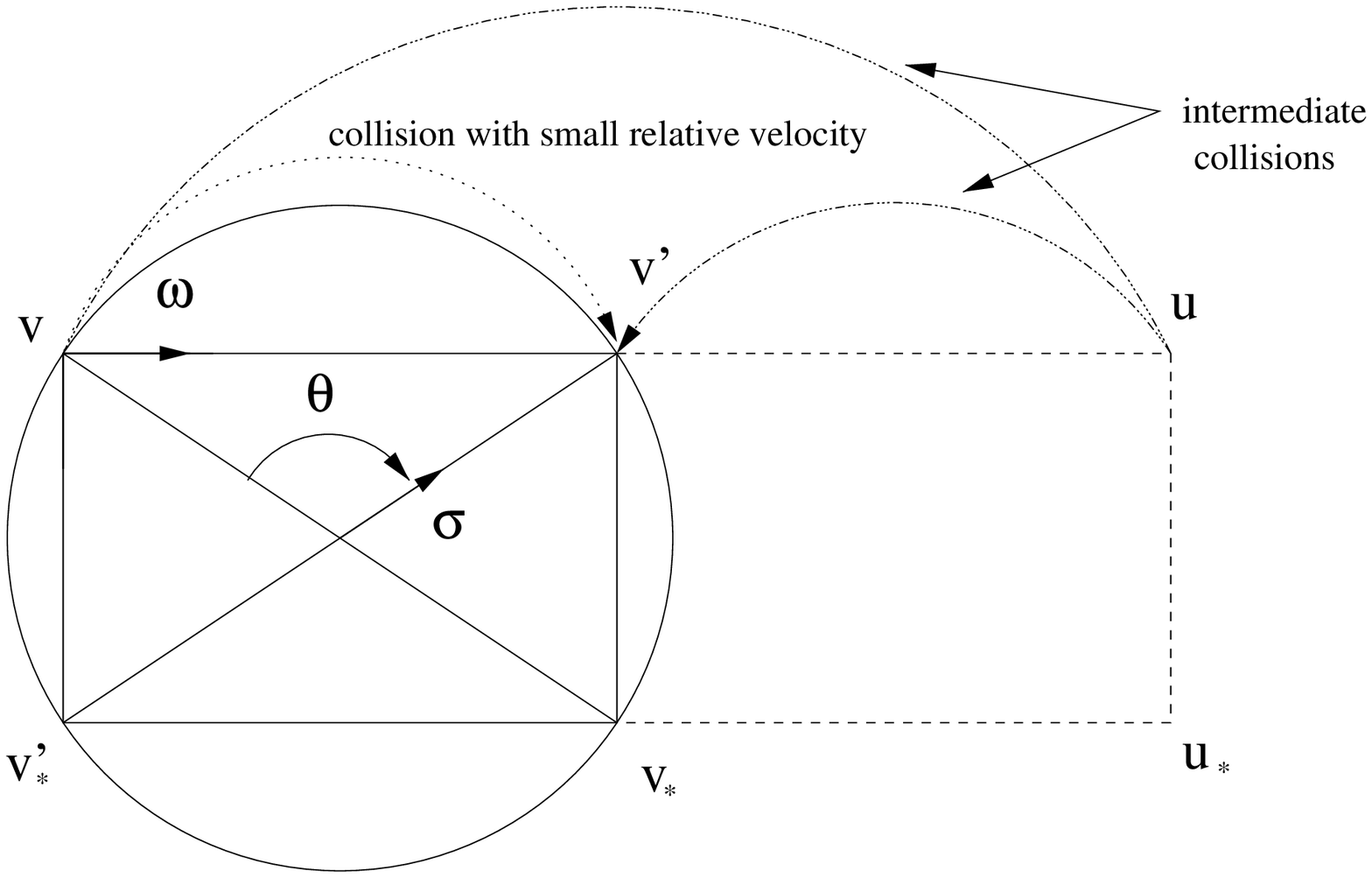}$$
 \caption{Introduction of an intermediate collision}\label{fig:2D}
 \end{figure}

Then, we shall use the following ``triangular'' inequality on the collision points:
 \begin{eqnarray*}
 \big[ \left( h(v) + h(v_*) \right) - \left( h(v') + h(v' _*) \right) \big]^2 
 &\le& 2 \big[ \left( h(v) + h(u_*) \right) - \left( h(u) + h(v_* ') \right) \big]^2 \\
 && + 2 \big[ \left( h(u) + h(v_*) \right) - \left( h(v') + h(u_*) \right) \big]^2.
 \end{eqnarray*}

Recall that $\int_{\R} e^{-r^2} \, dr = \sqrt{\pi}$. Let us add a third 
artificial integration variable $r_3$ on $\R$
 \begin{eqnarray*} 
 D^{{\ron Bo}} _1 (h) &\ge& \frac{1}{4\, \sqrt{\pi}} \int_{\ens{S}^{N-1}} \int_{V_1 \in \omega^\bot} 
 \int_{V_2 \in \omega^\bot} e^{-\left|V_1\right| ^2 -\left|V_2\right| ^2} 
 \int_{r_1 \in \R} \int_{r_2 \in \R} \int_{r_3 \in \R} 
 \Phi(|r_2 - r_1|) \, \tilde{b}(\theta_{1,2}) \\ 
 &&  k(r_1 \omega + V_1, r_2 \omega + V_2, r_2 \omega + V_1, r_1 \omega + V_2) \,
 e^{-r_1 ^2- r_2 ^2 - r_3 ^2} \, dr_1 \, dr_2 \, dr_3 \, dV_1 \, dV_2 \, d\omega. 
 \end{eqnarray*}
From now on, indexes of $\theta$ denote the points which are chosen to compute the angle. 
Now we rename $r_1, r_2, r_3$ first in $r_1, r_3, r_2$, secondly in $r_3, r_2, r_1$ 
and we take the mean of these two quantities. We get
 \begin{eqnarray*} 
 D^{{\ron Bo}} _1 (h) &\ge& \frac{1}{8 \sqrt{\pi}} \int_{\ens{S}^{N-1}} \int_{V_1 \in \omega^\bot} 
 \int_{V_2 \in \omega^\bot} e^{-|V_1|^2 - |V_2| ^2} \int_{r_1 \in \R}
 \int_{r_2 \in \R} \int_{r_3 \in \R} e^{-r_1 ^2-r_2 ^2 - r_3 ^2} \\ 
 && \Big[ \tilde{b}(\theta_{1,3}) \, \Phi(|r_3 - r_1|) \, 
 k(r_1 \omega + V_1, r_3 \omega + V_2, r_3 \omega + V_1, r_1 \omega + V_2) \\ 
 &&+ \tilde{b}(\theta_{2,3}) \, \Phi(|r_2 - r_3|) 
 \, k(r_3 \omega + V_1, r_2 \omega + V_2, r_2 \omega + V_1, r_3 \omega + V_2) \Big] 
 \, dr_1 \, dr_2 \, dr_3 \, dV_1 \, dV_2 \, d\omega. 
 \end{eqnarray*}
Then,
 \begin{eqnarray}\label{eq:zobi}
 D^{{\ron Bo}} _1 (h) &\ge& \frac{1}{8\, \sqrt{\pi}} 
 \int_{\ens{S}^{N-1}} \int_{V_1 \in \omega^\bot}
 \int_{V_2 \in \omega^\bot} e^{-|V_1|^2 - |V_2|^2} \int_{r_1 \in \R} 
 \int_{r_2 \in \R} \int_{r_3 \in \R}  \\ \nonumber
 && \min \big\{ \tilde{b}(\theta_{1,3}) \, \Phi(|r_3 - r_1|), 
 \tilde{b}(\theta_{2,3}) \, \Phi(|r_2 - r_3|) \big\} \\ \nonumber
 && \Big[ k(r_1 \omega + V_1, r_3 \omega + V_2, r_3 \omega + V_1, r_1 \omega + V_2) + \\ \nonumber 
 && k(r_3 \omega + V_1, r_2 \omega + V_2, r_2 \omega + V_1, r_3 \omega + V_2) \Big] \, 
 e^{-r_1 ^2- r_2 ^2 -r_3 ^2} \, dr_1 \, dr_2 \, dr_3 \, dV_1 \, dV_2 \, d\omega. 
 \end{eqnarray}
Now we use the following triangular inequality above-mentioned 
which means translated on $k$
 \begin{eqnarray*}\nonumber
 k(r_1 \omega + V_1, r_2 \omega + V_2, r_2 \omega + V_1, r_1 \omega + V_2) &\le&   
 2 \, k(r_1 \omega + V_1, r_3 \omega + V_2, r_3 \omega + V_1, r_1 \omega + V_2) \\ \nonumber
 && + 2 \, k(r_3 \omega + V_1, r_2 \omega + V_2, r_2 \omega + V_1, r_3 \omega + V_2). 
 \end{eqnarray*}
Plugging it in~\eqref{eq:zobi}, we get
 \begin{eqnarray*} \nonumber
 D^{{\ron Bo}} _1 (h) &\ge& \frac{1}{16 \, \sqrt{\pi}} 
 \int_{\ens{S}^{N-1}} \int_{V_1 \in \omega^\bot}
 \int_{V_2 \in \omega^\bot} e^{-|V_1|^2 - |V_2|^2} \int_{r_1 \in \R} 
 \int_{r_2 \in \R} \int_{r_3 \in \R}  \\ \nonumber
 && \min \big\{ \tilde{b}(\theta_{1,3}) \, \Phi(|r_3 - r_1|), 
 \tilde{b}(\theta_{2,3}) \, \Phi(|r_2 - r_3|) \big\} \\ \nonumber
 && k(r_1 \omega + V_1, r_2 \omega + V_2, r_2 \omega + V_1, r_1 \omega + V_2) \, 
 e^{-r_1 ^2- r_2 ^2 -r_3 ^2} \, dr_1 \, dr_2 \, dr_3 \, dV_1 \, dV_2 \, d\omega, 
 \end{eqnarray*}
which yields 
 \begin{eqnarray*} \nonumber
 D^{{\ron Bo}} _1 (h) &\ge& \frac{1}{16 \sqrt{\pi}} \int_{\ens{S}^{N-1}} \int_{V_1 \in \omega^\bot}
 \int_{V_2 \in \omega^\bot} e^{-|V_1|^2 - |V_2|^2} \int_{r_1 \in \R} 
 \int_{r_2 \in \R} \\ \nonumber
 && \left( \int_{r_3 \in \R} \min \big\{ \tilde{b}(\theta_{1,3}) \, \Phi(|r_3 - r_1|), 
 \tilde{b}(\theta_{2,3}) \, \Phi(|r_2 - r_3|) \big\} 
 \, e^{-r_3 ^2} \, dr_3 \right) \\
 && k(r_1 \omega + V_1, r_2 \omega + V_2, r_2 \omega + V_1, r_1 \omega + V_2) \, 
 e^{-r_1 ^2- r_2 ^2} \, dr_1 \, dr_2 \, dV_1 \, dV_2 \, d\omega. 
 \end{eqnarray*}
We now restrict the domain of integration for $r_3$ to the set  
 \[ \cal{D}_{r_1,r_2} = \big\{ r_3 \in \R \ | \ |r_3 - r_1| \ge |r_1 - r_2| \mbox{ and } 
 |r_2 - r_3| \ge |r_1 - r_2| \big\}. \]
Since $\tilde{b}$ is non-decreasing, 
and 
 \[ \cos \theta = \frac{|V_1 - V_2|^2 - |r_1 - r_2|^2}{|V_1 - V_2|^2 + |r_1 - r_2|^2} \]
which is non-increasing with respect to $|r_1 - r_2|$ when $V_1,V_2$ are kept frozen, 
it is easy to check 
that on this domain we have $\theta_{1,3} \ge \theta_{1,2}$ and 
$\theta_{2,3} \ge \theta_{1,2}$ and thus 
$\tilde{b}(\theta_{1,3}) \ge \tilde{b}(\theta_{1,2})$ and 
$\tilde{b}(\theta_{2,3}) \ge \tilde{b}(\theta_{1,2})$. Therefore we get 
 \begin{eqnarray*} \nonumber
 D^{{\ron Bo}} _1 (h) &\ge& \frac{1}{16 \sqrt{\pi}} \int_{\ens{S}^{N-1}} \int_{V_1 \in \omega^\bot}
 \int_{V_2 \in \omega^\bot} e^{-|V_1|^2 - |V_2|^2} \int_{r_1 \in \R} 
 \int_{r_2 \in \R} \\ \nonumber
 && \left( \int_{r_3 \in \cal{D}_{r_1,r_2}} \min \big\{\Phi(|r_3 - r_1|), \Phi(|r_2 - r_3|) \big\} 
 \, e^{-r_3 ^2} \, dr_3 \right) \tilde{b}(\theta_{1,2}) \\
 && k(r_1 \omega + V_1, r_2 \omega + V_2, r_2 \omega + V_1, r_1 \omega + V_2) \, 
 e^{-r_1 ^2- r_2 ^2} \, dr_1 \, dr_2 \, dV_1 \, dV_2 \, d\omega. 
 \end{eqnarray*}
Under assumption~\eqref{eq:hypPhi}, an easy computation leads to 
 \[ \left( \int_{r_3 \in \cal{D}_{r_1,r_2}} \min \big \{\Phi(|r_3 - r_1|),\Phi(|r_2 - r_3|) \big\} 
 \, e^{-|r_3|^2} \, dr_3 \right) \ge c_\Phi \, \sqrt{\pi} e^{- 4 R^2} >0 \]
as soon as $|r_1 - r_2| \le R$, i.e 
 \[ \left( \int_{r_3 \in \cal{D}_{r_1,r_2}} \min \big \{\Phi(|r_3 - r_1|),\Phi(|r_2 - r_3|) \big\} 
 \, e^{-|r_3|^2} \, dr_3 \right) \ge c_\Phi \, \sqrt{\pi} e^{-4R^2} \, 1_{|r_1 - r_2| \le R}. \]
By taking the mean of this estimate and the one obtained by replacing $\Phi$ by 
its bound from below $c_\Phi \, 1_{r \ge R}$, we deduce that
 \begin{eqnarray*}
 D^{{\ron Bo}} _1 (h) &\ge& \min \left( \frac{c_\Phi \, e^{-4R^2}}{8}, 
 \frac{c_\Phi}{2} \right) \, \frac{1}{4} \, 
 \int_{\ens{S}^{N-1}} \int_{V_1 \in \omega^\bot}
 \int_{V_2 \in \omega^\bot} e^{-|V_1|^2 - |V_2|^2} \int_{r_1 \in \R} 
 \int_{r_2 \in \R} \\ 
 && k(r_1 \omega + V_1, r_2 \omega + V_2, r_2 \omega + V_1, r_1 \omega + V_2) \, \tilde{b}(\theta) \, 
 e^{-r_1 ^2- r_2 ^2} \, dr_1 \, dr_2 \, dV_1 \, dV_2 \, d\omega.
 \end{eqnarray*}
If we now go back to the classical representation and simplify the minimum, we obtain
 \begin{eqnarray*} 
 D^{{\ron Bo}} _1 (h) &\ge& \left( \frac{c_\Phi \, e^{-4R^2}}{8} \right) \,   
 \frac{1}{4} \, \int_{\ens{S}^{N-1}} \int_{\R^N}\int_{\R^N} 
 M\, M_{*} \, k(v,v_*,v',v' _*) 
 \, d\sigma \, dv \, dv_* \\ 
 &=& \left( \frac{c_\Phi \, e^{-4R^2}}{8} \right) \,  D^{{\ron Bo}} _0 (h)
 \end{eqnarray*}
which concludes the proof of the lemma.
\end{proof}
The proof of Theorem~\ref{theo:Bolt} is a straightforward consequence of inequalities~\eqref{eq:ineg1} 
and~\eqref{eq:ineg2}.

\section{The Landau linearized operator} \label{sec:Landau}
\setcounter{equation}{0}

We now prove Theorem~\ref{theo:Land}. 
The idea here is to take the grazing collisions limit in some inequalities on 
the Boltzmann linearized operator obtained thanks to the geometrical method used 
in Section~\ref{sec:Bolt}. 
In fact, the most natural idea would have been to look for a geometrical property 
on the Landau linearized operator similar to the triangular inequality used for 
the Boltzmann linearized operator. 
But as collision circles reduce to lines in the grazing limit, 
the triangular inequality becomes trivial, and thus does not seem sufficient 
to apply the method of section~\ref{sec:Bolt}. It could be linked to the fact 
that in the grazing collisions limit one loses some information on the geometry of the collision. 
 
The problem that has to be tackled is to keep track of the angular collision kernel 
$b$. In fact we need it only for particular $b$, namely 
 \begin{equation}\label{eq:seqb}
 b_\var (\theta) = \frac{j_\var(\theta)}{\var^2 \sin^{N-2} \frac{\theta}{2}}
 \end{equation}
where $j_\var(\theta)=j(\theta/\var)/\var$ is a sequence of 
mollifiers (approximating $\delta_{\theta=0}$) with compact support in
$[0,\pi/2]$ and non-increasing on this interval. It is easy to see that $\tilde{b}_\var = 
2^{N-1} \sin^{N-2} \left( \frac{\theta}{2} \right) b_\var$ is also non-increasing on $[0,\pi]$. 
Following the same strategy as in Lemma~\ref{lem:phi} but keeping track 
of the angular part of the collision kernel, one obtains 
 \begin{lemma}\label{lem:phibis}
 Under the assumptions~\eqref{eq:hyptens}, \eqref{eq:hypPhi}, \eqref{eq:hypb}, 
 plus the assumption that $\tilde{b} = 2^{N-1} \sin^{N-2} \left( \frac{\theta}{2} \right) b$ 
 is non-increasing, one gets for all $h \in L^2(M)$ 
  \begin{equation}\label{eq:reducB:bis}
  D^{{\ron Bo}} _{b,\Phi} (h) \ge 
  \left( \frac{c_\Phi \, \beta_R}{8 \alpha_N} \right)\, D^{{\ron Bo}} _{b,1} (h)
  \end{equation}
 with 
  \[ \alpha_N = \int_{\R^{N-1}} e^{-|V|^2} \, dV, \ \ \ 
  \beta_R = \int_{\big\{V \in \R^{N-1} \ | \ |V| \ge 2 R\big\}} e^{-|V|^2} \, dV. \]
 Here $D^{{\ron Bo}} _{b,\Phi}$ stands for the entropy dissipation functional with $B = \Phi \, b$ 
 and $D^{{\ron Bo}} _{b,1}$ stands for the entropy dissipation functional with $B = b$.
 \end{lemma}
\begin{proof}[Proof of Lemma~\ref{lem:phibis}]
The geometrical idea of Lemma~\ref{lem:phi} can be applied to the variables $V_1,V_2$. 
Let us thus introduce two collisions points $u$ and $u_*$ (see figure~\ref{fig:2D:bis}) and replace
the collision ``$(v,v_*)$ gives $(v',v'_*)$'' by the two collisions
``$(v,u_*)$  gives $(v',u)$'' and  ``$(u,v_*)$ gives $(u_*,v'_*)$''.
 \begin{figure}[h]
 \epsfysize=6cm
 $$\epsfbox{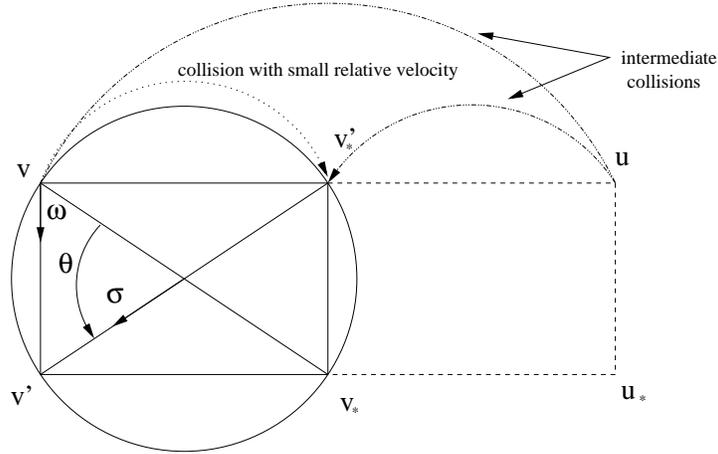}$$
 \caption{Introduction of an intermediate collision}\label{fig:2D:bis}
 \end{figure} 

Then, we shall use the following ``triangular'' inequality on the collision points:
 \begin{eqnarray}
 \big[ \left( h(v) + h(v_*) \right) - \left( h(v') + h(v' _*) \right) \big]^2 
 &\le& 2 \big[ \left( h(v) + h(u_*) \right) - \left( h(v') + h(u) \right) \big]^2 \nonumber\\
 && + 2 \big[ \left( h(u) + h(v_*) \right) - \left( h(u_*) + h(v' _*) \right) \big]^2 \nonumber. 
 \end{eqnarray}
Now we introduce an artificial third variable $V_3$ on $\omega^\bot$. Let us denote 
 \begin{equation*} 
 \alpha_N = \int_{\R^{N-1}} e^{-|V|^2} \, dV
 \end{equation*}
By inverting either $V_1$ and $V_3$ or $V_2$ and $V_3$, taking the mean, and using 
the ``triangular'' inequality above-mentioned we get 
 \begin{eqnarray*} \nonumber
 D^{{\ron Bo}} _{b,\Phi} (h) &\ge& 
 \frac{1}{16 \alpha_N} \,\int_{\ens{S}^{N-1}} \int_{r_1 \in \R} \int_{r_2 \in \R} e^{-r_1 ^2-r_2 ^2} 
 \int_{V_1 \in \omega^\bot} \int_{V_2 \in \omega^\bot} \\ \nonumber
 && \left( \int_{V_3 \in \omega^\bot} \min \left\{ \tilde{b} (\theta_{1,3}) \, \Phi(|V_3 - V_1|),
 \tilde{b} (\theta_{2,3}) \, \Phi(|V_2 - V_3|)\right\} 
 \, e^{-|V_3|^2} \, dV_3 \right) \\
 && k(r_1 \omega + V_1, r_2 \omega + V_2, r_2 \omega + V_1, r_1 \omega + V_2) \, e^{-|V_1|^2-|V_2|^2} 
 dr_1 \, dr_2 \, dV_1 \, dV_2 \, d\omega.
 \end{eqnarray*}
Let us now restrict the integration along $V_3$ to the domain 
 \[\cal{D}_{V_1,V_2} = \big\{ V_3 \ | \ |V_3 - V_1| \ge |V_1 - V_2| \mbox{ and } 
 |V_2 - V_3| \ge |V_1 - V_2| \big\}. \] 
Then since the expression 
 \[ \cos \theta = \frac{|V_1 - V_2|^2 - |r_1 - r_2|^2}{|V_1 - V_2|^2 + |r_1 - r_2|^2} \] 
is non-decreasing according to $|V_1 - V_2|$ when $r_1,r_2$ are kept frozen, and 
$\tilde{b}$ is non-increasing, we get $\theta_{1,3} \le \theta_{1,2}$ and 
$\theta_{2,3} \le \theta_{1,2}$ (see figure~\ref{fig:2D:bis}) and 
so $\tilde{b} ( \theta_{1,3}) \ge \tilde{b} ( \theta_{1,2})$ and 
$\tilde{b} ( \theta_{2,3}) \ge \tilde{b} ( \theta_{1,2})$. Consequently 
  \begin{eqnarray*} \nonumber
 D^{{\ron Bo}} _{b,\Phi} (h) &\ge& 
 \frac{1}{16 \alpha_N} \,\int_{\ens{S}^{N-1}} \int_{r_1 \in \R} \int_{r_2 \in \R} dr_2 e^{-r_1 ^2-r_2 ^2} 
 \int_{V_1 \in \omega^\bot} \int_{V_2 \in \omega^\bot} \\ \nonumber
 && \left( \int_{V_3 \in \cal{D}_{V_1,V_2}} \min \left\{ \Phi(|V_3 - V_1|), \Phi(|V_2 - V_3|)\right\} 
 \,  e^{-|V_3|^2} \, dV_3 \right) \, \tilde{b} (\theta_{1,2}) \\
 && k(r_1 \omega + V_1, r_2 \omega + V_2, r_2 \omega + V_1, 
 r_1 \omega + V_2) \, e^{-|V_1|^2-|V_2|^2} 
 dr_1 \, dr_2 \, dV_1 \, dV_2 \, d\omega.
 \end{eqnarray*}
Under assumption~\eqref{eq:hypPhi}, an easy computation leads to 
 \begin{multline*} 
 \left( \int_{V_3 \in \cal{D}_{V_1,V_2}} \min \big \{\Phi(|V_3 - V_1|),\Phi(|V_2 - V_3|) \big\} 
 \, e^{-|V_3|^2} \, dV_3 \right) \\
 \ge  c_\Phi \, \int_{\big\{V \in \R^{N-1} \ | \ |V| \ge 2 R \big\}} e^{-|V|^2} \, dV = c_\Phi \, \beta_R >0
 \end{multline*}
as soon as $|V_1 - V_2| \le R$, i.e 
 \[ \left( \int_{V_3 \in \cal{D}_{V_1,V_2}} \min \big \{\Phi(|V_3 - V_1|),\Phi(|V_2 - V_3|) \big\} 
 \, e^{-|V_3|^2} \, dV_3 \right) \ge c_\Phi \, \beta_R \, 1_{|V_1 - V_2| \le R} \, . \]
Taking the mean of this estimate and 
the one obtained by the trivial lower bound $\Phi(r) \ge c_\Phi 1_{\{r \ge R\}}$, 
we get in the end 
\begin{equation*}
 D^{{\ron Bo}} _{b,\Phi} (h) \ge \min \left( \frac{c_\Phi \, \beta_R}{8 \alpha_N} ,  
 \frac{c_\Phi}{2} \right) \frac{1}{4} \int_{\ens{S}^{N-1}} \int_{\R^N} \int_{\R^N} 
 b (\theta) \, M \, M_{*} \, \left[ h' _* + h' - h - h_* \right]^2 
 \, d\sigma \, dv_{*} \, dv,\nonumber 
 \end{equation*}
which yields
 \begin{equation*}
 D^{{\ron Bo}} _{b,\Phi} (h) 
 \ge \left( \frac{c_\Phi \, \beta_R}{8 \alpha_N} \right) \, D^{{\ron Bo}} _{b,1} (h)
 \end{equation*}
and concludes the proof of the lemma. 
\end{proof}

We now have to take the grazing collisions limit in the entropy dissipation functional to 
prove inequality~\eqref{eq:reduc:Landau} of 
Theorem~\ref{theo:Land} (this limit is essentially well-known, see for instance~\cite{Desv:asBE:92}). 
 \begin{lemma}\label{lem:limite} 
 Let us consider $\, b_\var$ as defined in~\eqref{eq:seqb} and $\Phi$ satisfying 
 assumption~\eqref{eq:hypPhi}. Then for a given $h \in L^2 (M)$, 
  \begin{equation*}
  D^{{\ron Bo}} _{b_\var,\Phi} (h) \xrightarrow[\var \to 0]{} c_{N,j} \, D^{{\ron La}} _\Phi (h)
  \end{equation*}
 where
  \[ c_{N,j} = \frac{2^{N-5} \,|\ens{S}^{N-2}|}{N-1} \, 
     \left( \int_o ^\pi j(\chi) \, \chi^2 \, d\chi \right) \]
 depends only on the dimension $N$ and the mollifier $j$. 
 $D^{{\ron Bo}} _{b_\var,\Phi}$ stands for the Boltzmann entropy dissipation functional with 
 $B = \Phi \, b_\var$, and $D^{{\ron La}} _\Phi$ stands for the Landau entropy dissipation 
 functional with collision kernel $\Phi$. 
 \end{lemma}
\begin{proof}[Proof of Lemma~\ref{lem:limite}] 
The idea of the proof is to expand the expression for small $\var$ and is very similar 
to what is done in~\cite{Desv:asBE:92}. Let us write the angular vector $\sigma$
 \begin{equation*}
 \sigma=\frac{v-v_*}{|v-v_*|}\cos(\theta)+\mathbf{n} \,\sin(\theta),
 \end{equation*}
where $\mathbf{n}$ is a unit vector in  $(v-v_*)^\bot$. Therefore, we shall write
 \begin{eqnarray}
 D^{{\ron Bo}} _{b_\var,\Phi} (h) & = &\frac{1}{4} \int_{\R^N}\int_{\R^N} \Phi(|v-v_*|) \, M\, M_{*}
 \int_{\ens{S}^{N-2}((v-v_*)^\bot)} \int_{\theta=0}^{\pi}  \, b_\var(\theta) \, \nonumber\\
 &&\Bigg[ h \left( v-\frac{v-v_*}{2}(1-\cos(\theta)) +
 \frac{|v-v_*|}{2}\,\mathbf{n} \,\sin(\theta) \right)\nonumber\\
 &&+h\left(v_*+\frac{v-v_*}{2}(1-\cos(\theta))-\frac{|v-v_*|}{2}\,\mathbf{n}\,\sin(\theta)\right)\nonumber\\
 &&-h(v)-h(v_*)\Bigg]^2 \,\sin^{N-2} \theta \, d\theta\, d\mathbf{n} \, dv \, dv_*,\nonumber
 \end{eqnarray}
where $\ens{S}^{N-2}((v-v_*)^\bot)$ denotes the unit sphere in $(v-v_*)^\bot$. 
Let us now focus on the integral on $\theta$
 \begin{eqnarray*}
 && \int_{\theta=0}^{\pi}  \, b_\var(\theta) \, 
 \Bigg[ h \left( v-\frac{v-v_*}{2}(1-\cos(\theta)) +
 \frac{|v-v_*|}{2}\,\mathbf{n} \,\sin(\theta) \right) \\
 &&+h\left(v_*+\frac{v-v_*}{2}(1-\cos(\theta))-\frac{|v-v_*|}{2}\,\mathbf{n}\,\sin(\theta)\right) 
 -h(v)-h(v_*)\Bigg]^2 \,\sin^{N-2} \theta \, d\theta,
 \end{eqnarray*}
and make the change of variables $\chi=\theta/\var$. We get 
 \begin{eqnarray*}
 && \int_{\chi=0}^{\pi} \frac{\sin^{N-2}(\var\,\chi)}{\sin^{N-2}(\frac{\var\,\chi}{2})} 
 \, \frac{j(\chi)}{\var^2} \Bigg[ h \left( v-\frac{v-v_*}{2}(1-\cos(\var\,\chi)) +
 \frac{|v-v_*|}{2}\, \mathbf{n}\,\sin(\var\,\chi) \right) \\
 &&+h\left(v_*+\frac{v-v_*}{2}(1-\cos(\var\,\chi))-\frac{|v-v_*|}{2}\, \mathbf{n} \,\sin(\var\,\chi)\right) 
 -h(v)-h(v_*)\Bigg]^2 \, d\chi
 \end{eqnarray*}
i.e for small $\var$,
 \begin{equation*}
 \int_{\chi=0}^{\pi} (2^{N-2} + O(\var)) \, \frac{j(\chi)}{\var^2} 
 \, \left( \frac{|v-v_*|}{2} \right)^2 \, 
 \Big[\var\,\chi \, \mathbf{n} \, \cdot 
 \left( \nabla_v h(v)- \nabla_{v_*} h(v_*) \right)+O(\var^2 \chi^2 )\Big]^2 \, d\chi, 
 \end{equation*}
which writes 
 \begin{multline*}
 |v-v_*|^2 \, \int_{\chi=0}^{\pi} 2^{N-4} \, j(\chi) \, \chi^2 \, \Big[ \mathbf{n} \, \cdot 
 \left( \nabla_v h(v)- \nabla_{v_*} h(v_*) \right)\Big]^2 \, d\chi + O(\var) \\
 = 2^{N-4} \, \left( \int_o ^\pi j(\chi) \, \chi^2 \, d\chi \right) \, |v-v_*|^2 \, 
 \Big[ \mathbf{n} \, \cdot \left( \nabla_v h(v)- \nabla_{v_*} h(v_*) \right)\Big]^2 + O(\var).
 \end{multline*}
As the unit vector $\mathbf{n}$ is orthogonal to $\frac{v-v_*}{|v-v_*|}$, we can 
introduce here the orthogonal projection onto $(v-v_*)^\bot$
 \begin{multline*}
 D^{{\ron Bo}} _{b_\var,\Phi} (h) = \frac{2^{N-4}}{4} \, 
 \left( \int_o ^\pi j(\chi) \, \chi^2 \, d\chi \right) \, \int_{\R^N} \int_{\R^N}
 \int_{\ens{S}^{N-2}((v-v_*)^\bot)} \Phi(|v-v_*|) \, |v-v_*|^2  \\ 
 \Big[ \mathbf{n} \cdot \Pi_{(v-v_*)^\bot} 
 \big( \nabla_v h(v)- \nabla_{v_*} h(v_*) \big) \Big]^2 M \, M_* \, d\mathbf{n} \, dv \, dv_* + O(\var). 
 \end{multline*}
It is straightforward to see that 
 \begin{equation*} 
 \int_{\ens{S}^{N-2}((v-v_*)^\bot)} \big( \mathbf{n} \cdot u \big)^2 \, d\mathbf{n} 
 = \zeta_N \left\|  u \right\|^2
 \end{equation*} 
with, for any $\mathbf{u} \in \ens{S}^{N-2}$
  \[ \zeta_N = \int_{\ens{S}^{N-2}} (\mathbf{u} \cdot \mathbf{n})^2 \, d\mathbf{n} 
         = \frac{|\ens{S}^{N-2}|}{N-1}. \]
Thus we get in the end
 \begin{eqnarray*}
 D^{{\ron Bo}} _{b_\var,\Phi} (h) & = & \frac{|\ens{S}^{N-2}| \, 2^{N-4}}{4 (N-1)} \, 
 \left( \int_o ^\pi j(\chi) \, \chi^2 \, d\chi \right) \\ 
 && \int_{\R^N}\int_{\R^N}
 |v-v_*|^2 \, \Phi(|v-v_*|) \, M\, M_{*} 
 \left\| \nabla_v h(v)- \nabla_{v_*} h(v_*) \right\|^2 \, dv_{*}\, dv + O(\var) \\
 &=& \frac{|\ens{S}^{N-2}| \, 2^{N-5}}{N-1} \, \left( \int_o ^\pi j(\chi) \, \chi^2 \, d\chi \right) 
 \, D^{{\ron La}} _\Phi (h) + O(\var). 
 \end{eqnarray*}
This concludes the proof of Lemma~\ref{lem:limite}.
\end{proof}

Coming back to the proof of Theorem~\ref{theo:Land}, 
we first prove~\eqref{eq:reduc:Landau}: we write down 
inequality~\eqref{eq:reducB:bis} on $D^{{\ron Bo}} _{\Phi,b_\var}$ 
since $\tilde{b}_\var$ is non-increasing, and we apply Lemma~\ref{lem:limite} 
on each term, which gives 
 \begin{eqnarray*}
 D^{{\ron La}} (h) \ge C^{{\ron La}} _\Phi \, D^{{\ron La}}_{0} (h),
 \end{eqnarray*}
where 
 \[ C^{{\ron La}} _\Phi = \left( \frac{c_\Phi \, \beta_R}{8 \alpha_N} \right). \]
Inequality~\eqref{eq:tsL} follows immediatly.

It remains to prove the lower 
bound~\eqref{eq:vpL} on the first non-zero eigenvalue of the Landau linearized 
operator for Maxwellian molecules in dimension $3$. Let us denote by 
$\lambda^{{\ron Bo}} _{0,b_\var}$ the first non-zero eigenvalue for the Boltzmann 
linearized operator with $B=b_\var$: for all $h \in L^2(M)$ orthogonal in $L^2(M)$ 
to $1, v, |v|^2$,  
 \[ D^{{\ron Bo}} _{b_\var} (h) \ge |\lambda^{{\ron Bo}} _{0,b_\var}| 
 \, \|h\|^2 _{L^2(M)}. \]
We apply Lemma~\ref{lem:limite} to this inequality which leads to 
 \[ D^{{\ron La}} _0 (h) \ge 
    \frac{ \lim_{\var \to 0} |\lambda^{{\ron Bo}} _{0,b_\var}|}{c_{3,j}} \,\|h\|^2 _{L^2(M)} \] 
for all $h \in L^2(M)$ orthogonal in $L^2(M)$ to $1, v, |v|^2$.
An explicit formula for $|\lambda^{{\ron Bo}} _{0,b_\var}|$ is given in~\cite{Boby:maxw:88} 
 \[ |\lambda^{{\ron Bo}} _{0,b_\var}| = 
    \pi \, \int_0 ^\pi \sin^3 (\theta) \, b_\var (\theta) \, d\theta \]
and thus 
 \[ \lim_{\var \to 0} |\lambda^{{\ron Bo}} _{0,b_\var}| = 
     2\pi \, \left( \int_0 ^\pi j(\chi) \, \chi^2 \, d\chi \right) \]
which concludes the proof.
\bigskip

\noindent
{\bf{Acknowledgment}} : Support by the European network HYKE, funded by the EC as
contract HPRN-CT-2002-00282, is acknowledged. 

\medskip

\bibliographystyle{acm}

\bibliography{LBE}

\begin{flushleft} \signcb \end{flushleft}
\vspace*{-55mm} \begin{flushright} \signcm \end{flushright}

\end{document}